\documentclass[12pt,oneside,draft]{amsart}
\usepackage[T1]{fontenc}
\usepackage[latin1]{inputenc}
\usepackage{geometry}
\usepackage{amssymb}

\makeatletter
 \theoremstyle{plain}
\newtheorem{thm}{Theorem}[section]
  \theoremstyle{plain}
  \newtheorem{lem}[thm]{Lemma}
  \theoremstyle{plain}
  \newtheorem{cor}[thm]{Corollary}
 \theoremstyle{definition}
  \newtheorem{example}[thm]{Example}
  \theoremstyle{remark}
  
  \theoremstyle{definition}
  \newtheorem{defn}[thm]{Definition}
  \theoremstyle{plain}
  \newtheorem{prop}[thm]{Proposition}

\makeatletter
\def\Na{\mathcal N}
\def\Ka{\mathcal K}
\def\Ja{\mathcal J}
\def\LL{\mathbb L}

\def\Ca{\mathcal C}
\def\R{\RR}
\def\Fa{\mathcal F}
\def\a{\alpha}
\def\tr{\mbox{tr}}
\def\mun{\hat\mu^N}

\def\cxm{\mathbb C\langle X_1,\cdots,X_m\rangle}
\def\D{\Delta}
\def\G{\Gamma}
\def\RR{\mathbb R}
\def\d{\delta}
\def\e{\epsilon}

\def\ra{\rightarrow}

\makeatother

\usepackage[english]{babel}
\makeatother
\begin{document}
%\selectlanguage{english}

\title{On  Classical Analogues of Free Entropy Dimension.}

\author{A. Guionnet{*} and D. Shlyakhtenko{*}{*}}

\begin{abstract}
We define  a classical probability analogue of Voiculescu's free entropy
dimension that we shall call
the classical probability  entropy
dimension of a probability measure on $\mathbb{R}^n$. We show that
the classical probability entropy
dimension of a measure  is related with diverse other notions
of dimension. First, it can be viewed as a kind of 
fractal dimension. 
Second, if one extends Bochner's inequalities 
to a  measure by requiring that microstates
around this measure 
 asymptotically satisfy the classical  Bochner's inequalities,
then we show that  the classical probability entropy
dimension controls 
the rate of increase
of optimal constants in Bochner's inequality for a measure regularized
by convolution with the Gaussian law as the regularization is removed.
We introduce a free analogue of the Bochner inequality and study the
related free entropy dimension quantity. We show that it is
greater or equal to  the non-microstates free entropy dimension.
\end{abstract}

\thanks{{*}Ecole Normale Sup\'erieure de Lyon, Unit\'e de Math\'ematiques
pures et appliqu\'ees, UMR 5669, 46 All\'ee d'Italie, 69364 Lyon
Cedex 07, France and the Miller institute for
Basic Research in Science, University of California Berkeley. E-mail: aguionnet@umpa.ens-lyon.fr}

\thanks{{*}{*}Department of Mathematics, UCLA, Los Angeles, CA 90095. E-mail:
shlyakht@math.ucla.edu. Research supported by NSF grants DMS-0355226
and DMS-0555680.}

\maketitle

\section{Introduction.}

In \cite{Vo2}, using his notion of free entropy $\chi$, Voiculescu
introduced the free entropy dimension of a non-commutative law. If
$X_{1},\ldots,X_{n}\in(M,\tau)$ are self-adjoint non-commutative
random variables in a tracial $W^{*}$-probability space, then\[
\delta(X_{1},\ldots,X_{n})=n+\limsup_{t\to0}\frac{\chi(X_{1}^{t},\ldots,X_{n}^{t})}{|\log t|},\]
 where $X_{j}^{t}=X_{j}+tS_{j}$ and $S_{1},\ldots,S_{n}$ form a
free semicircular family, free from $X_{1},\ldots,X_{n}$. Voiculescu's
motivation was to introduce a kind of asymptotic Minkowski content
of matricial microstate spaces associated to the joint law of $X_{1},\ldots,X_{n}$.
Indeed, for a variation of the definition of free entropy dimension,
K. Jung has proved a formula that involves asymptotic packing numbers
\cite{Jung}. Moreover, he proved (again, for a version of the definition
above), that one obtains the same number whether one uses semicircular
perturbations or some other perturbation $X_{j}^{t}=X_{j}+tY_{j}$,
where $Y_{1},\ldots,Y_{n}$ are some $n$-tuple, free from $X_{1},\ldots,X_{n}$
and having finite free entropy.

The free entropy dimension is a remarkable quantity, with unexpected
connections to other branches of mathematics. For example, if $X_{1},\ldots,X_{n}$
generate the group algebra of a discrete group $\Gamma$, $\delta(X_{1},\ldots,X_{n})$
is related by an inequality to the $L^{2}$-Betti numbers of the group
$\Gamma$ (this is based on a number of results, see \cite{CSh,MSh}).
Unfortunately, the exact values of free entropy dimension are known
in only a few cases. For example, in the case of a single variable
$X$ with law given by a probability measure $\mu$ on $\mathbb{R}$,
$\delta(\mu)=1-\sum_{t\in\mathbb{R}}\mu(\{ t\})^{2}$.

One of the most important questions surrounding $\delta$ is the question
of its invariance under various functional calculi. It is hoped that
$\delta(X_{1},\ldots,X_{n})=\delta(Y_{1},\ldots,Y_{m})$ if $X_{1},\ldots,X_{n}$
and $Y_{1},\ldots,Y_{m}$ generate the same von Neumann algebra (i.e.,
are {}``non-commutative measurable functions of each other''). However,
the question is open even if it is asked for continuous functions
(that is, assuming that the $C^{*}$-algebras generated by $X_{1},\ldots,X_{n}$
and $Y_{1},\ldots,Y_{m}$ are the same). What is known, for a version
of the definition of free entropy dimension, is that its value is
preserved under algebraic changes of generators. Solving these problems
would be of great interest to von Neumann algebra theory.

In the first part of the present paper, we turn to look at the classical
analogue of free entropy dimension. Given a probability measure $\mu$
on $\mathbb{R}^{n}$ (which can be though of as the law of $n$ real
random variables $X_{1},\ldots,X_{n}$), we consider the measure $\mu_{t}=\mu*\nu_{t}$,
where $\nu_{t}$ is the Gaussian law\[
\nu_{t}(\prod dx_{j})=\frac{1}{(2\pi t^{2})^{n/2}}\exp(-\frac{1}{2t^{2}}\sum x_{j}^{2})\prod dx_{j}.\]
 Thus $\mu_{t}$ is the law of $X_{1}^{t},\ldots,X_{n}^{t}$ with
$X_{j}^{t}=X_{j}+tG_{j}$, and $G_{1},\ldots,G_{n}$ independent Gaussian
random variables, independent from $X_{1},\ldots,X_{n}$. We then
set\[
\delta_{c}(\mu)=n-\liminf_{t\to0}\frac{H(\mu_{t})}{|\log t|},\]
 where for a non negative Lebesgue absolutely-continuous measure $p(x)dx$,
\[
H(p(x)dx)=\int p(x)\log p(x)dx.\]
 (The change of sign here is due to the fact that $H(\mu_{t})\in(-\infty,+\infty]$
behaves as the analogue of $-\chi$).

The main result of this paper relates $\delta_{c}(\mu)$ with a kind
of average fractal dimension of the measure $\mu$. In particular,
we prove that $\delta_{c}(\mu)$ remains the same if $\mu$ is replaced
by a push-forward by a Lipschitz function. However, the value of $\delta_{c}(\mu)$
may change if we push forward $\mu$ by a continuous or measurable
function.

We also prove a number of technical properties of $\delta_{c}$. Among
the ones of independent interest is the fact that (in the case that
$\limsup$ in its definition is a limit) $\delta_{c}$ is affine:
$\delta_{c}(\sum\alpha_{j}\mu_{j})=\sum\alpha_{j}\delta_{c}(\mu_{j})$
in the case that $\mu_{j}$ are probability measures and $\alpha_{j}\geq0$,
$\sum\alpha_{j}=1$.

The second part of the paper relates the rate of increase of optimal
constants in 
an ad hoc notion of
Bochner's inequality for measures with entropy dimension. We say that
a probability measure $\mu$ satisfies Bochner's inequality with constants
$(n, K(n))\in ({\mathbb R}^+)^2$ if for all smooth $f$,
\begin{equation}\label{tyu}
\mu(\Gamma_{2}(f,f))\geq\frac{1}{n}\mu((\Delta f)^{2})-K(n)\mu(\Gamma(f,f)),
\end{equation}
where $\Gamma(f,f)$ and $\Gamma_{2}(f,f)$ are the carr\'e du champ
and carr\'e du champ it\'er\'e, respectively. Intuitively, one
should think of $n$ as the dimension of the support of $\mu$ and
$K$ as an estimate for the smallest eigenvalue of the Ricci curvature
of the support in the
sense that if $\mu=\d_x$,
we recover the classical Bochner inequality
at the point $x$, with $n$ the
dimension of the manifold
where $x$ lives and $-K(n)$
a lower bound on the Ricci curvature (cf. e.g. \cite{BQ,BL}). 
The definition is actually obtained by considering
the microstates $\Gamma_
N(\mu,\e):=\{x_1,\cdots,x_N\in\mathbb R^N:d(N^{-1}\sum_{i=1}^N \d_{x_i},
\mu)<\e\}$, viewing it as a submanifold of $\R^N$
with some dimension $[nN]$ and Ricci curvature bounded below by
 $-K(n)$. Letting then $N$ going to infinity
gives \eqref{tyu}. 
We now replace $\mu$ with
$\mu_{\varepsilon}=\mu*\nu_{\varepsilon}$ and study the 
 functions
$\varepsilon\mapsto K(n,\varepsilon)\geq0$ such that $\mu_{\sqrt{\e}}$
satisfies  Bochner's inequality with constants
$(n,K(n,\e))$.  We then set\[
\delta^{\square}=1-\inf(\liminf_{\varepsilon\to0}\frac{\int_{\varepsilon}^{1}K(n,y)dy}{|\log\varepsilon|}+1)n,\]
where the $\inf$ is taken over all $n\ge 0$
and  functions $K(n,\varepsilon)$
for which \eqref{tyu} holds. We prove that with this definition,
$\delta^{\square}=\delta_{c}$.

In the third and final part of the paper, we study the free non-commutative
analogue of the inequality \eqref{tyu} and the related free
entropy dimension quantity, which we show to be less than or equal
to the non-microstates free entropy dimension.

\section{Equivalent definitions of $\delta_{c}$.}

The main result of this section is that one can replace in the definition
of $\delta_{c}(\mu)$ the convolution with the Gaussian measure by
convolution with dilations of any other probability measure $\nu$
that has finite entropy. We first consider some properties of $\delta_{c}$,
which are of independent interest. Throughout this section, it will
be convenient to assume that $\nu$ is a finite positive measure,
but to drop the assumption that its total mass is $1$. We will also
denote by $D_{t}:\mathbb{R}\to\mathbb{R}$ the dilation map $x\mapsto tx$.
For simplicity of notation, we give all statements and proofs for
a measure on $\mathbb{R}$. However, these go through unaltered for
measures on $\mathbb{R}^{n}$. Also, all $\liminf$ could be replaced
by $\limsup$ if one would prefer to define $\delta_{c}$ with a $\limsup$.

\begin{lem}
\label{lem:bound}(a) Let $\nu$ be a Lebesgue absolutely continuous
finite measure on $\mathbb{R}$, $\nu_{t}=D_{t}^{*}(\nu)$ (where
$D_{t}$ is the map $x\mapsto tx$ is a dilation). Then for any probability
measure $\mu$ and any constant $\alpha>0$ we have\[
\liminf_{t\to0}\frac{H(\alpha\mu_{t})}{\log t}=\alpha\liminf_{t\to0}\frac{H(\mu_{t})}{|\log t|}\]
 \\
 (b) Let $\nu$ be a non negative Lebesgue absolutely continuous measure
for which $\nu(\mathbb{R})=\delta<\infty$. Let $\mu$ be a probability
measure on $\mathbb{R}$ and denote $\nu_{t}=D_{t}^{*}(\nu)$ and
$\mu_{t}=\mu*\nu_{t}$. Let $p_{t}(x)$ be the density of $\mu_{t}$.

If \begin{equation}
\int\log(1+|x|)d\nu(x)<\infty,\,\mbox{ and }\,\int\log(1+|x|)d\mu(x)<\infty,\label{tight}\end{equation}
 then \[
0\leq\liminf_{t\to0}\frac{H(\mu_{t})}{|\log t|}.\]

On the other hand, if $H(\nu)<\infty$, then \[
\liminf_{t\to0}\frac{H(\mu_{t})}{|\log t|}\leq\limsup_{t\to0}\frac{H(\mu_{t})}{|\log t|}\leq\delta.\]

(Here and below $H(q(x)dx)=\int q(x)\log q(x)dx$ for any non-negative
measurable function $q$, even if $q(x)dx$ is not a probability measure). 
\end{lem}
\begin{proof}
(a) follows from the formula $H(\alpha\mu)=\alpha H(\mu)+\mu(\mathbb{R})\log\alpha$
and the fact that $\mu_{t}(\mathbb{R})=\nu(\mathbb{R})$ is independent
of $t\in\mathbb{R}$.

(b)We may assume without loss of generality that $\delta=1$ by a
rescaling up to using (a).

For the first inequality, recall that for any probability measure
$\nu$, any non negative function $f$, Jensen's inequality implies
that \[
\int f(x)\log f(x)d\nu(x)\ge\int f(x)d\nu(x)\log\left(\int f(x)d\nu(x)\right).\]
 Therefore, if we let $\nu(dx)=p(x)dx$ be a probability measure absolutely
continuous with respect to the Lebesgue measure, we can write \[
H(f(x)dx)=\int\frac{f(x)}{p(x)}\log\frac{f(x)}{p(x)}p(x)dx+\int\log p(x)f(x)dx\ge\int\log p(x)f(x)dx\]
 if $\int f(x)dx=1$. We can for instance take $p(x)=\frac{1}{2(1+|x|)^{2}}$
to obtain the lower bound \[
H(f(x)dx)\ge-\log2-2\int\log(1+|x|)f(x)dx\]
 for all $f\ge0$ so that $\int f(x)dx=1$.

Now, since $\nu$ is absolutely continuous with respect to Lebesgue
measure, so is the measure $\mu_{t}(dx)=f_{t}(x)dx$. Applying the
above to $f_{t}$, we deduce \begin{eqnarray*}
H(\mu_{t}) & \ge & -\log2-2\int\log(1+|x|)d\mu_{t}(x)\\
 & \ge & -\log2-2\int\log(1+|x|)(1+t|y|)d\mu(x)d\nu(y)\\
 & \ge & -\log2-2\int\log(1+|x|)d\mu(x)-2\int\log(1+|y|)d\nu(y)\end{eqnarray*}
 where the last bound holds for $t\le1$. Hence, when \eqref{tight}
is satisfied, $H(\mu_{t})$ is bounded below independently of $t\le1$,
which gives the desired lower bound.

We next prove the upper bound. By the entropy power inequality (see
e.g. \cite{stam}), we have that\begin{eqnarray*}
\exp(-2H(\mu_{t})) & \geq & \exp(-2H(\mu))+\exp(-2H(\nu_{t}))\\
 & \geq & \exp(-2H(\nu_{t}))\\
 & = & \exp(-2H(\nu)+2\log t).\end{eqnarray*}
 Thus\[
H(\mu_{t})\leq H(\nu)-\log t\]
 so that\[
\limsup_{t\to0}\frac{H(\mu_{t})}{|\log t|}\leq\limsup_{t\to0}\frac{H(\nu)-\log t}{|\log t|}=1\]
 as claimed. 
\end{proof}
\begin{lem}
\label{lem:affimemu}Let $n\in\mathbb{N}$ and $\mu=\sum_{i=1}^{n}\mu_{i}$
for some non negative measures $(\mu_{i},1\le i\le n)$ so that $\mu_{i}(\mathbb{R})=a_{i}>0$,
$\sum_{i=1}^{n}a_{i}=1$. %Assume that $\mu$ and $\nu$ have finite variance.
Let $\nu$ be a probability measure on $\mathbb{R}$ so that $H(\nu)<\infty$.
Then \begin{equation}
\liminf_{t\to0}\frac{H(\mu*\nu_{t})}{|\log t|}=\liminf_{t\to0}\frac{1}{|\log t|}\sum a_{i}H(a_{i}^{-1}\mu_{i}*\nu_{t}).\label{gogo}\end{equation}

\end{lem}
Note that since $H(\nu)$ is assumed finite, $H(a_{i}^{-1}\mu_{i}*\nu_{t})\le|\log t|$
by the previous Lemma and so the sum in the right hand side of \eqref{gogo}
is well defined.

\begin{proof}
Since $\nu$ is absolutely continuous with respect to Lebesgue measure
with density $p$, so is $\mu*\nu_{t}$ and \[
p_{\mu}(x)=\frac{d\mu*\nu_{t}}{dx}(x)=\frac{1}{t}\int p\left(\frac{x-y}{t}\right)d\mu(y).\]
 We assume first that $n=2$ and denote in short $p_{i}(x)=a_{i}^{-1}p_{\mu_{i}}(x)$
for $i=1,2$, so that $\int p_{i}(x)dx=1$. Then the density of $\mu*\nu_{t}$
is given by $\sum a_{i}p_{i}(x)$ and hence\begin{eqnarray*}
H(\mu*\nu_{t}) & = & \int\sum_{i}a_{i}p_{i}(x)\log\sum_{j}a_{j}p_{j}(x)\ dx=\sum_{i}a_{i}\int p_{i}(x)\log(\sum_{j}a_{j}p_{j}(x))dx.\end{eqnarray*}
 As a consequence, \[
H(\mu*\nu_{t})-\sum a_{i}H(a_{i}^{-1}\mu_{i}*\nu_{t})=\sum_{i}a_{i}\int p_{i}(x)\log\left(\frac{\sum_{j}a_{j}p_{j}(x)}{a_{i}p_{i}(x)}\right)dx+\sum_{i=1}^{2}a_{i}\log a_{i}.\]
 Then for each $i=1,2$\[
\sum_{j}a_{j}p_{j}(x)/a_{i}p_{i}(x)=1+\frac{a_{j}p_{j}(x)}{a_{i}p_{i}(x)}\]
 where in the last term $i,j\in\{1,2\}$ and $i\neq j$.

Since for $y\geq0$, $0\leq\log(1+y)\leq y$ and since $p_{j}(x)$,
$p_{i}(x)\geq0$, we conclude that\[
0\leq\log\left(1+\frac{a_{j}p_{j}(x)}{a_{i}p_{i}(x)}\right)\leq\frac{a_{j}p_{j}(x)}{a_{i}p_{i}(x)}.\]
 Hence\[
0\leq H(\mu*\nu_{t})-\sum_{i=1}^{2}a_{i}H(a_{i}^{-1}\mu_{i}*\nu_{t})\leq\sum_{j}\int a_{j}p_{j}(x)dx+\sum a_{i}\log a_{i}\leq1+\sum a_{i}\log a_{i}.\]
 If $\mu=\sum_{i=1}^{n}\mu_{i}$ for $n>2$, we first apply the above
bound with $\mu_{1}'=\mu_{1},\mu_{2}'=\sum_{i=2}^{n}\mu_{i}$ and
$a_{1}'=a_{1}$, $a_{2}'=\sum_{i=2}^{n}a_{i}$, and then proceed by
induction, replacing $\mu$ by $(\sum_{i=2}^{n}a_{i})^{-1}\sum_{i=2}^{n}\mu_{i}$.
We get in this way

\[
0\leq H(\mu*\nu_{t})-\sum_{i=1}^{n}a_{i}H(a_{i}^{-1}\mu_{i}*\nu_{t})\leq n-1+\sum_{i=1}^{n}a_{i}\log a_{i}.\]
 Thus\[
\lim_{t\to0}\frac{H(\mu*\nu_{t})-\sum_{i=1}^{n}a_{i}H(a_{i}^{-1}\mu_{i}*\nu_{t})}{|\log t|}=0,\]
 which implies the claim. 
\end{proof}
We have as an immediate corollary a somewhat surprising property of
$\delta_{c}$:

\begin{cor}
Assume that $\mu_{j}$ are probability measures for which $\limsup$
in the definition of $\delta_{c}$ is a limit. Then the map $\mu\mapsto\delta_{c}(\mu)$
is affine: if $\alpha_{j}\geq0$, $\sum\alpha_{j}=1$, then $\delta_{c}(\sum\alpha_{j}\mu_{j})=\sum\alpha_{j}\delta_{c}(\mu_{j})$. 
\end{cor}
Note that this property is very particular to the commutative case.
Indeed, recall that the formula for the free entropy dimension of
a single self-adjoint variable with law $\mu$ can be equivalently written as\[
\delta(\mu)=1-\sum_{t\in\mathbb{R}}\mu\times\mu(\{(t,t)\})\]
 so that $\delta(\mu)$ is \emph{quadratic} in $\mu$. By the Cauchy-Schwartz
inequality, one has $\delta(\sum_{i=1}^{n}a_{i}\mu_{i})\ge\sum_{i=1}^{n}a_{i}\delta(\mu_{i})$
but equality can hold only if for all $t\in\mathbb{R}$, $\mu_{i}(\{ t\})$
does not depend on $i\in\{1,\cdots,n\}$.

\begin{lem}
\label{lem:affinenu}Let for $n\in\mathbb{N}$, $\nu=\sum_{i=1}^{n}\nu^{(i)}$
so that $\nu^{(i)}(\mathbb{R})=a_{i}$. Assume that $H(a_{i}^{-1}\nu^{(i)})$
is finite for all $i$. Then\[
\liminf_{t\to0}\frac{H(\mu*\nu_{t})}{|\log t|}=\liminf_{t\to0}\frac{1}{|\log t|}\sum a_{i}H(a_{i}^{-1}\mu*\nu_{t}^{(i)}).\]

\end{lem}
\begin{proof}
The proof is very similar to that of Lemma \ref{lem:affimemu} and
we first assume $n=2$. We let $\nu_{t}^{(i)}=D_{t}^{*}\nu^{(i)}$
where $D_{t}:\mathbb{R}\to\mathbb{R}$ is the map $D_{t}(x)=tx$.
We have:\[
\mu*\nu_{t}=\sum_{i}\mu*\nu_{t}^{(i)}=\sum_{i}a_{i}(a_{i}^{-1}\mu*\nu_{t}^{(i)}).\]
 Thus if we set\[
p_{i}(x)=d(a_{i}^{-1}\mu*\nu_{t}^{(i)})/dx\]
 then the density of $\mu*\nu_{t}$ is given by $\sum a_{i}p_{i}(x)$
and hence\begin{eqnarray*}
H(\mu*\nu_{t}) & = & \int\sum_{i}a_{i}p_{i}(x)\log\sum_{j}a_{j}p_{j}(x)\ dx=\sum_{i}a_{i}\int p_{i}(x)\log(\sum_{j}a_{j}p_{j}(x))dx.\end{eqnarray*}
 Hence, we deduce as in the proof of Lemma \ref{lem:affimemu} that
\[
0\leq H(\mu*\nu_{t})-\sum a_{i}H(p_{i}(x)dx)\leq\sum_{j}\int a_{j}p_{j}(x)dx+\sum a_{j}\log a_{j}\leq1+\sum a_{j}\log a_{j}.\]
 Thus\[
\lim_{t\to0}\frac{H(\mu*\nu_{t})-\sum a_{i}H(p_{i}(x)dx)}{|\log t|}=0,\]
 which implies the claim. 
\end{proof}
\begin{cor}
\label{corr:cutoff}Given $\nu(dx)=f(x)dx$, with $\nu(\mathbb{R})=1$
and $H(\nu)<\infty$, set $\nu_{t}=D_{t}^{*}(\nu)$ where $D_{t}:\mathbb{R}\to\mathbb{R}$,
given by $D_{t}(x)=tx$. Let $\mu$ be a probability measure on $\mathbb{R}$
%of finite variance
. Then given $\varepsilon>0$ there exists $M$ sufficiently large
so that if we denote by $\nu^{M}$ the measure $\nu([-M,M])^{-1}\nu|_{[-M,M]}$,
$\nu_{t}^{M}=D_{t}(\nu^{M})$ and by $\mu_{M}$ the measure $\mu_{M}=\mu[-M,M]^{-1}\mu|_{[-M,M]}$,
then\[
\left|\liminf_{t\to0}\frac{H(\mu*\nu_{t})}{|\log t|}-\liminf_{t\to0}\frac{H(\mu_{M}*\nu_{t}^{M})}{|\log t|}\right|<\varepsilon.\]

\end{cor}
\begin{proof}
This follows from first decomposing $\mu$ as $\mu|_{[-M,M]}+\mu|_{[-M,M]^{c}}$,
so that Lemma \ref{lem:affimemu} shows that \begin{multline*}
\left|\liminf_{t\to0}\frac{H(\mu*\nu_{t})}{|\log t|}-\mu([-M,M])\liminf_{t\to0}\frac{H(\mu_{M}*\nu_{t})}{|\log t|}\right|\\
\le\mu([-M,M]^{c})\limsup_{t\to0}\frac{H(\mu([-M,M]^{c})^{-1}\mu|_{[-M,M]^{c}}*\nu_{t})}{|\log t|}\le\mu([-M,M]^{c})\end{multline*}
 where the last inequality is due to Lemma \ref{lem:bound}.(b) since
$\nu(\mathbb{R})=1$.

We next decompose $\nu$ as $\nu|_{[-M,M]}+\nu|_{[-M,M]^{c}}$ and
apply Lemma \ref{lem:affinenu}. Since $H(\nu)$ is finite, also $H(\nu|_{[-M,M]})$
and $H(\nu|_{[-M,M]^{c}})$ are finite and so \begin{multline*}
\left|\liminf_{t\to0}\frac{H(\mu_{M}*\nu_{t})}{|\log t|}-\nu([-M,M])\liminf_{t\to0}\frac{H(\mu_{M}*\nu_{t}^{M})}{|\log t|}\right|\\
\le\nu([-M,M]^{c})\limsup_{t\to0}\frac{H(\mu|_{[-M,M]}*D_{t}^{*}(\nu_{M})}{|\log t|}
\le\nu([-M,M]^{c})\end{multline*}
 again by Lemma \ref{lem:bound}(b). Since \[
\left|(\mu([-M,M])\nu([-M,M])-1)\liminf_{t\to0}\frac{H(\mu_{M}*\nu_{t}^{M})}{|\log t|}\right|\le\mu([-M,M]^{c})+\nu([-M,M]^{c})\]
 the proof is complete if we take $M$ big enough so that $2(\mu([-M,M]^{c})+\nu([-M,M]^{c}))\le\varepsilon$. 
\end{proof}
\begin{lem}
\label{lemma:replaceconst}Assume that $\nu(dx)=f(x)dx$ with $\mbox{supp}f=E$
a bounded subset of $\mathbb{R}$, and that for some constant $C>\varepsilon>0$,
$|f-C|<\varepsilon$ on $E$. Let $\nu'(dx)=C\chi_{E}dx$ and set
$\nu_{t}=D_{t}^{*}(\nu)$, $\nu'_{t}=D_{t}^{*}(\nu')$. Assume furthermore
that the support of $\mu$ is a bounded subset of $\mathbb{R}$. Then
\[
\left|\liminf_{t\to0}\frac{H(\mu*\nu_{t}')}{|\log t|}-\liminf_{t\to0}\frac{H(\mu*\nu_{t})}{|\log t|}\right|\leq\varepsilon\lambda(E).\]

\end{lem}
\begin{proof}
Recall that \begin{gather*}
p_{t}(x):=\frac{d\mu*\nu_{t}}{dx}(x)=\int f(t^{-1}(x-y))\frac{1}{t}d\mu(y)\\
p_{t}'(x):=\frac{d\mu*\nu_{t}'}{dx}(x)=C\int\chi_{E}(t^{-1}(x-y))\frac{1}{t}d\mu(y).\end{gather*}

Using the fact that $\mu$ is a probability measure, we have:\[
\left|p_{t}(x)-p_{t}'(x)\right|\leq\int\varepsilon\chi_{E}(t^{-1}(x-y))\frac{1}{t}d\mu(y)=\frac{\varepsilon}{C}p_{t}'(x).\]
 In particular, we have that\[
\left|\frac{p_{t}(x)}{p_{t}'(x)}-1\right|\leq C^{-1}\varepsilon.\]
 Thus\[
\int p_{t}(x)\log p_{t}(x)dx=\int p_{t}(x)\log p_{t}'(x)dx-\int p_{t}(x)\log\frac{p_{t}(x)}{p'_{t}(x)}dx,\]
 implies \[
|\int p_{t}(x)\log p_{t}(x)dx-\int p_{t}(x)\log p_{t}'(x)dx|\leq\max|\log(1\pm C^{-1}\varepsilon)|=f(C^{-1}\varepsilon),\]
 with $f(C^{-1}\varepsilon)\to0$ as $C^{-1}\varepsilon\to0$. Hence\begin{eqnarray*}
|H(\mu*\nu_{t}')-H(\mu*\nu_{t})| & \leq & |\int(p_{t}(x)-p_{t}'(x))\log p_{t}'(x)dx|+f(C^{-1}\varepsilon)\\
 & \leq & \frac{\varepsilon}{C}\int p_{t}'(x)|\log p_{t}'(x)|dx+f(C^{-1}\varepsilon).\end{eqnarray*}
 It follows that \begin{equation}
\left|\liminf_{t\to0}\frac{H(\mu*\nu_{t}')}{|\log t|}-\liminf_{t\to0}\frac{H(\mu*\nu_{t})}{|\log t|}\right|\leq\frac{\varepsilon}{C}\limsup_{t\to0}\frac{\int p_{t}'(x)|\log p_{t}'(x)|dx}{|\log t|}.\label{popo}\end{equation}

Now, let $A_{t}=\{ x:0<p_{t}'(x)\leq1\}\subset tE+\textrm{supp}\mu$.
Then $\log p_{t}'(x)>0$ for $x\notin A_{t}$ and $\log p_{t}'(x)\leq0$
for $x\in A_{t}$. Therefore,\[
\int p_{t}'(x)|\log p_{t}'(x)|dx=\int p_{t}'(x)\log p_{t}'(x)dx-2\int_{A_{t}}p_{t}'(x)\log p_{t}'(x)dx.\]
 Since for $y\in[0,1]$, the function $y\log y$ is bounded from below
by $-e^{-1}$ and from above by $0$, we get that for $x\in A_{t}$,
$0\le-p_{t}'(x)\log p_{t}'(x)\leq e^{-1}$. Since $A_{t}\subset E+\textrm{supp}\mu$
for $t\le1$, the Lebesgue measure $\lambda(A_{t})$ is bounded uniformly
in $t$. Thus, we find that \[
\liminf_{t\to0}\frac{\int p_{t}'(x)|\log p_{t}'(x)|dx}{|\log t|}=\liminf_{t\to0}\frac{\int p_{t}'(x)\log p_{t}'(x)}{|\log t|}=\liminf_{t\to0}\frac{H(\mu*\nu_{t}')}{|\log t|}.\]
 But since $H(\nu')={C\lambda(E)}\log C$ is finite, we can use Lemma
\ref{lem:bound} to conclude that the right hand side above is bounded
by $C\lambda(E)$, the mass of $\nu$. Hence, we have proved with
\eqref{popo} that \[
\left|\liminf_{t\to0}\frac{H(\mu*\nu_{t}')}{|\log t|}-\liminf_{t\to0}\frac{H(\mu*\nu_{t})}{|\log t|}\right|\leq\varepsilon\lambda(E).\]

\end{proof}
\begin{thm}
\label{thm:anynu}Let $\nu$ be an arbitrary probability measure with
$H(\nu)$ finite. Assume that $\mu$ is a probability measure, and
assume that $\mu$ and $\nu$ satisfy \eqref{tight}.% of finite variance. 
Then if we denote by $D_{t}^{*}$ the push-forward of a measure by
the dilation $x\mapsto tx$, we have that\[
\liminf_{t\to0}\frac{H(\mu*D_{t}^{*}(\nu))}{|\log t|}=\liminf_{t\to0}\frac{H(\mu*D_{t}^{*}(\chi_{[0,1]}))}{|\log t|}.\]
 In particular, the limit is independent of the measure $\nu$. 
\end{thm}
\begin{proof}
Fix $\varepsilon>0$. By Corollary \ref{corr:cutoff}, we may assume,
without changing $\liminf_{t\to0}\frac{H(\mu*D_{t}^{*}(\nu))}{|\log t|}$
by more than $\varepsilon/2$, that $\mu$ and $\nu$ are supported
on bounded sets. In particular, $\nu$ is Lebesgue absolutely continuous
with density $q(x)\in L^{1}(\mathbb{R})$ with $E=\textrm{supp}q$
a subset of finite Lebesgue measure. Given $\varepsilon>0$ we may
find a subset $E_{0}\subset\mathbb{R}$ and a constant $M$ so that
$q(x)<M$ on $E_{0}$ and $\nu(E_{0})^{-1}\leq1-\varepsilon/8$. By
Corollary \ref{corr:cutoff} we may replace $\nu$ by $\nu(E_{0})^{-1}\nu|_{E_{0}}$
without affecting the value of $\liminf_{t\to0}\frac{H(\mu*D_{t}^{*}(\nu))}{|\log t|}$
by more than $\varepsilon/4$. Next, since the density $p(x)$ of
$\nu$ is now a bounded function on the support of $\nu$, we may
find a finite collection of disjoint subsets $E_{j}\subset E_{0}$
and constants $C_{j}$ with the property that on each $E_{j}$, $|p_{j}-C_{j}|<\varepsilon/\lambda(E)8$
and that $C_{j}$ is the average value of $f$ on $E_{j}$ (in particular,
$\sum C_{j}\lambda(E_{j})=\int f(x)dx=1$). According to Lemma \ref{lemma:replaceconst}
we may replace on each $E_{j}$ $\nu|_{E_{j}}$ with $\chi_{E_{j}}$
at a penalty of at most $\varepsilon\lambda(E_{j})/8$. Hence we may
replace $\nu$ with the probability measure $\sum C_{j}\chi E_{j}$
at a penalty of at most $(\varepsilon\lambda(E)/8)\cdot\sum\lambda(E_{j})\leq\varepsilon/8$.
By Lemma \ref{lem:affinenu} it follows that\[
\liminf_{t\to0}\frac{H(\mu*\nu_{t})}{|\log t|}=\liminf_{t\to0}\sum\frac{H(\mu*D_{t}^{*}(C_{j}\chi_{E_{j}}))}{|\log t|}.\]

Finally, by Lebesgue almost everywhere differentiability theorem,
we may find, for each $E_{j}$ disjoint intervals $I_{1}^{(j)},\dots,I_{k_{j}}^{(j)}$
of rational length with the property that $E_{j}$ and $\cup_{k}I_{k}^{(j)}$
differ by at most $\lambda(E_{j})\cdot\varepsilon/8$. Applying once
again Lemma \ref{lem:affimemu} and Lemma \ref{lem:affinenu}, we
conclude that we may assume at a further penalty of $\varepsilon/8$
that $ $$\nu=\sum K_{r}\chi_{E_{r}}$ where $E_{r}$ are a finite
collection of intervals. Up to subdivision, we may assume that all
the $E_{r}$ have the same Lebesgue measure (or length). We conclude
that\[
\liminf_{t\to0}\frac{H(\mu*\nu_{t})}{|\log t|}=\liminf_{t\to0}\sum K_{r}\frac{H(\mu*D_{t}^{*}(\chi_{E_{r}}))}{|\log t|}+o(\varepsilon),\]
 where $K_{r}$ is a family of non negative real numbers so that $\sum K_{r}\lambda(E_{r})=1$
and $E_{r}$ are intervals. 

Since $H(q(x)dx)=H(q(x-y)dx)$, we may replace any interval $E_{r}$
in the previous formula by a shifted interval $E_{j}+k_{j}$ for any
constant $k_{j}$. Hence, since all the $E_{r}$ have the same length,
$H(\mu*D_{t}^{*}(\chi_{E_{r}}))$ does not depend on $r$ and so we
have \[
\liminf_{t\to0}\frac{H(\mu*\nu_{t})}{|\log t|}=\liminf_{t\to0}\frac{1}{\lambda(E_{1})}\frac{H(\mu*D_{t}^{*}(\chi_{E_{1}}))}{|\log t|}+o(\varepsilon),\]
 here $E_{1}$ is an interval with right hand point at the origin.
Note that $E_{1}$ could be chosen as small as wished and so letting
$\varepsilon$ going to zero we have 
\[
\liminf_{t\to0}\frac{H(\mu*\nu_{t})}{|\log t|}=\lim_{a\downarrow0}\liminf_{t\to0}\frac{H(\mu*D_{t}^{*}(\chi_{[0,a]}))}{a|\log t|}.\]
 This shows in particular that $\liminf_{t\to0}\frac{H(\mu*\nu_{t})}{|\log t|}$
does not depend on the probability measure $\nu$ with finite entropy
and so we also have \[
\liminf_{t\to0}\frac{H(\mu*\nu_{t})}{|\log t|}=\liminf_{t\to0}\frac{H(\mu*D_{t}^{*}(\chi_{[0,1]}))}{|\log t|}.\]
 
\end{proof}

\section{$\delta_{c}$ and fractal dimension.}

If $\mu$ is a probability measure on $\mathbb{R}$, one can consider
the (lower) point wise dimension of $\mu$:\[
f^{\mu}(x)=\liminf_{t\to0}\frac{\mu[x-t,x+t]}{\log t}.\]
This function quantifies the logarithmic rate of growth of the measures
of $t$-balls around $x$ and hence is a kind of local fractal dimension
of $\mu$. For example, certain Cantor-Lebesgue measures\[
\mu=\frac{1}{2}(\delta_{-1}+\delta_{1})*\frac{1}{2}(\delta_{\lambda}+\delta_{-\lambda})*\frac{1}{2}(\delta_{\lambda^{2}}+\delta_{-\lambda^{2}})*\cdots,\qquad0<\lambda<1/2,\]
 satisfy $f^{\mu}=\alpha=-\log_{2}\lambda$ on the Cantor set supporting
$\mu$ and $f^{\mu}=0$ outside of it. We show that $\delta_{c}$
is very close to the average value (computed with respect to $\mu$)
of the function $f^{\mu}$, apart from the question of exchanging
integration against $\mu$ and the limit $\liminf_{t\to0}$.

\begin{thm}
\label{thm:fractal}Let $\mu$ be a probability measure on $\mathbb{R}$,
and let\[
d_{t}(x)=\frac{-\log\mu[x-t/2,x+t/2]}{|\log t|}.\]
 Then\[
\delta_{c}(\mu)=\limsup_{t\to0}\int d_{t}(y)d\mu(y).\]

\end{thm}
\begin{proof}
By Theorem \ref{thm:anynu} we may write\[
\delta_{c}(\mu)=1-\liminf_{t\to0}\frac{H(\mu*\nu_{t})}{|\log t|},\]
 where $\nu_{t}=D_{t}^{*}\chi_{[-1/2,1/2]}=\frac{1}{t}\chi_{[-t/2,t/2]}$.
Let $p_{t}(x)$ be the density of $\mu_{t}$:\[
p_{t}(x)=(\mu*D_{t}^{*}\chi_{[-1/2,1/2]})(x)=\frac{1}{t}\mu([x-t/2,x+t/2]).\]

Now,\begin{eqnarray*}
H(\mu_{t}) & = & \int p_{t}(x)\log p_{t}(x)dx\\
 & = & \iint\frac{1}{t}\chi_{[-t/2,t/2]}(x-y)d\mu(y)\log p_{t}(x)dx\\
 & = & \iint\frac{1}{t}\chi_{[-t/2,t/2]}(x)\log p_{t}(x+y)dxd\mu(y).\end{eqnarray*}
 Since $p_{t}(x+y)=\frac{1}{t}\mu[x+y-t/2,x+y+t/2]$ and $[y+x-t/2,y+x+t/2]\subset[y-t,y+t]$
as long as $-t/2\leq x\leq t/2$, we find that for $|x|\leq t/2$,
$p_{t}(x+y)\leq\frac{1}{t}\mu[y-t,y+t]$. Thus\begin{eqnarray*}
H(\mu_{t}) & \leq & \iint\frac{1}{t}\chi_{[-t/2,t/2]}(x)\log\frac{1}{t}\mu[y-t,y+t]\ dxd\mu(y)\\
 & = & \int\frac{1}{t}\chi_{[-t/2,t/2]}(x)dx\ \int\log\frac{1}{t}\mu[y-t,y+t]d\mu(y)\\
 & = & \int\log\frac{1}{2t}\mu[y-t,y+t]d\mu(y)+\int\log2d\mu(y)=\int\log\frac{1}{2t}\mu[y-t,y+t]+\log2\end{eqnarray*}
 (since $\mu$ is a probability measure). It follows that\[
\liminf_{t\to0}\frac{H(\mu_{t})}{|\log t|}\leq\liminf_{t\to0}\frac{\int\log\frac{1}{t}\mu[y-t/2,y+t/2]d\mu(y)}{|\log t|}=\liminf_{t\to0}\frac{\int\log p_{t}(y)d\mu(y)}{|\log t|}.\]

Let now $\delta>0$ and set $C=1+\delta$. Let $\nu'=\chi_{[-C/2,C/2]}$,
$\nu''=\nu'-\chi_{[-1/2,1/2]}$. Let $\mu_{t}'=\mu*D_{t}^{*}(\nu')$,
$\mu_{t}''=\mu*D_{t}^{*}(\nu'')$. Thus $\mu_{t}'=\mu_{t}+\mu_{t}''$.
Let $p'_{t}(x)$, $p_{t}''(x)$ be the densities of $\mu_{t}'$ and
$\mu_{t}''$, respectively. Then we have:\begin{eqnarray*}
\int p_{t}(x)\log p_{t}'(x)dx-\int p_{t}(x)\log p_{t}(x)dx & = & \int p_{t}(x)\log\frac{p'_{t}(x)}{p_{t}(x)}dx\\
 & = & \int p_{t}(x)\log\frac{p_{t}(x)+p_{t}''(x)}{p_{t}(x)}dx\\
 & = & \int p_{t}(x)\log(1+p_{t}''(x)/p_{t}(x))dx.\end{eqnarray*}
 Since $0\leq\log(1+z)\leq z$ for $z\geq0$, we conclude that\begin{eqnarray*}
0 & \leq & \int p_{t}(x)\log(1+p_{t}''(x)/p_{t}(x))dx\\
 & \leq & \int p_{t}(x)p_{t}''(x)/p_{t}(x)\ dx\\
 & = & \int p_{t}''(x)dx=\mu_{t}''(\mathbb{R})=\delta.\end{eqnarray*}
 It follows that \begin{equation}
\left|\int p_{t}(x)\log p_{t}'(x)dx-\int p_{t}(x)\log p_{t}(x)dx\right|\leq\delta.\label{ghgh}\end{equation}
 Now, $p_{t}'(x)=\frac{1}{t}\mu[x-Ct/2,x+Ct/2]$. If $|x|<t/2$, then
$[y-\delta t/2,y+\delta t/2]\subset[y+x-Ct/2,y+x+Ct/2]$. Thus $p_{t}'(x+y)\geq\frac{1}{t}\mu[y-\delta t/2,y+\delta t/2]$
as long as $|x|<t/2$. It follows that\begin{eqnarray}
\int p_{t}(x)\log p_{t}'(x)dx & = & \iint\frac{1}{t}\chi_{[-t/2,t/2]}(x)\log p_{t}'(x+y)\ d\mu(y)dx\nonumber \\
 & \geq & \iint\frac{1}{t}\chi_{[-t/2,t/2]}(x)\log\frac{1}{t}\mu[y-\delta t/2,y+\delta t/2]\ d\mu(y)dx\nonumber \\
 & = & \int\log\frac{1}{t}\mu[y-\delta t/2,y+\delta t/2]\ d\mu(y)\nonumber \\
 & = & \int\log\frac{1}{\delta t}\mu[y-\delta t/2,y+\delta t/2]\ d\mu(y)+\log\delta.\label{lklk}\end{eqnarray}
 Thus, first by \eqref{ghgh} and then \eqref{lklk} we obtain \begin{multline*}
\liminf_{t\to0}\frac{H(\mu_{t})}{|\log t|}=\liminf_{t\to0}\frac{\int p_{t}(x)\log p_{t}'(x)dx}{|\log t|}\\
\geq\liminf_{t\to0}\frac{\int\log\frac{1}{\delta t}\mu[y-\delta t/2,y+\delta t/2]\ d\mu(y)}{|\log t|}=\liminf_{t\to0}\frac{\int\log p_{t}(y)d\mu(y)}{|\log t|}\end{multline*}
 where we finally made the change of variable $t'=\delta t$. Combining
this with the previous estimate proves that\begin{eqnarray*}
\delta_{c}(\mu) & = & 1-\liminf_{t\to0}\frac{\int\log t^{-1}\mu[x-t/2,x+t/2]\ d\mu(x)}{|\log t|}\\
 & = & 1-\liminf_{t\to0}\int\left[\frac{\log\mu[x-t/2,x+t/2]}{|\log t|}+\frac{-\log t}{|\log t|}\right]\ d\mu(x)\\
 & = & \limsup_{t\to0}\int d_{t}(x)d\mu(x).\end{eqnarray*}

\end{proof}
\begin{cor}
Assume that $\mu$ is a probability measure, which is dimension regular;
i.e., there exists some $\mu$-measurable function $\alpha(x)$ and
strictly positive constants $C$, $c$, and $t_{0}$ so that for any
$x$ in the support of $\mu$ and all $0<t<t_{0}$ one has\begin{equation}
ct^{\alpha(x)}\leq\mu[x-t/2,x+t/2]\leq Ct^{\alpha(x)}.\label{eq:dim-regular}\end{equation}
 Then $\delta_{c}(\mu)=\int\alpha(x)d\mu(x)$. 
\end{cor}
Note that in all the previous results, we could have change the $\liminf$
into a $\limsup$ and vice versa. Under the hypotheses of the Corollary
we would thus obtain \[
\delta_{c}(\mu)=1-\liminf_{t\to0}\frac{H(\mu_{t})}{|\log t|}=1-\limsup_{t\to0}\frac{H(\mu_{t})}{|\log t|}=\int\alpha(x)d\mu(x).\]

\begin{proof}
We find that\[
d_{t}(x)=-\frac{\log\mu[x-t/2,x+t/2]}{|\log t|}\]
 satisfies the inequalities\[
\frac{\alpha(x)\log t+\log c}{|\log t|}\leq-d_{t}(x)\leq\frac{\alpha(x)\log t+\log C}{|\log t|},\]
 so that for $t<1$,\[
\alpha(x)-\frac{\log c}{|\log t|}\geq d_{t}(x)\geq\alpha(x)-\frac{\log C}{|\log t|}.\]
 Integrating these inequalities against $d\mu(x)$, passing to the
limit as $t\to0$ and using Theorem \ref{thm:fractal}, we obtain
that $\delta_{c}(\mu)=\int\alpha(x)d\mu(x)$. 
\end{proof}
\begin{example}
(i) Let $0<\alpha<1$ and let $\mu_{\alpha}$ be the Cantor-Lebesgue
measure given by\[
\mu_{\alpha}=\frac{1}{2}(\delta_{-1}+\delta_{1})*\frac{1}{2}(\delta_{\lambda}+\delta_{-\lambda})*\frac{1}{2}(\delta_{\lambda^{2}}+\delta_{-\lambda^{2}})*\cdots\qquad\lambda=2^{-\alpha}\]
 Then $\mu_{\alpha}$ satisfies \eqref{eq:dim-regular} with $\alpha(x)=\alpha$
for all $x$ in the support of $\mu_{\alpha}$. Thus $\delta_{c}(\mu_{\alpha})=\alpha$.\\
 (ii) Let $\mu=\delta_{0}$ be a delta measure at $0$. Then \eqref{eq:dim-regular}
is satisfied with $\alpha=0$ on the support of $\mu$. Hence $\delta_{c}(\mu)=0$.\\
 (iii) Let $\mu$ be Lebesgue absolutely continuous with density $p(x)$.
Then $\mu=\mu_{M}+\mu_{M}^{\perp}$ where $\mu_{M}=\mu|_{\{ x:p(x)\leq M\}}$.
Furthermore, $\mu_{M}^{\perp}(\mathbb{R})\to0$ as $M\to\infty$.
Thus by Lemma \eqref{lem:bound}, $\lim_{M\to\infty}\delta_{c}(\mu_{M}^{\perp})=0$
and hence  $\delta_{c}(\mu)=\lim_{M\to\infty}\delta_{c}(\mu_{M})+\delta_{c}(\mu_{M}^{\perp})=\lim_{M\to\infty}\delta_{c}(\mu_{M})$.
Since $H(\mu_{M})<\infty$ and by
the entropy power inequality $H(\mu_{M}*\nu)\leq H(\mu_{M})$
for any $\nu$, we find that 
 $\delta_{c}(\mu_{M})=1$. Thus $\delta_{c}(\mu)=1$. 
\end{example}
It is curious to note that one has a classical analogue of the connection
between free entropy dimension and group cohomology. In the classical
case, the $L^{2}$ Betti numbers are replaced with ordinary Betti
numbers and the statement greatly trivializes.

Let $\Gamma$ be a discrete abelian group, and let $\hat{\Gamma}$ be the 
its Pontrjagin dual 
$\hat\Gamma = \operatorname{Hom}(\Gamma,\{z\in\mathbb{C}:|z|=1\})$.
Then $\hat\Gamma$ is compact, and each $\gamma\in\Gamma$ can be identified
with a bounded function on $\hat\Gamma$ by $\gamma(\phi)=\phi(\gamma)$,
$\phi\in\hat\Gamma$.  Let $H^1(\Gamma,\mathbb{C})$ denote the group cohomology
of $\Gamma$ with coefficients in $\mathbb{C}$ (viewed as a trivial 
$\Gamma$-module).

\begin{thm}
Let $\Gamma$ be a finitely generated discrete abelian group with
generators $\gamma_{1},\ldots,\gamma_{n}$. Identify $\mathbb{C}\Gamma\subset L^{\infty}(\hat{\Gamma},\mu)$,
where $\mu$ is a Haar measure of $\hat{\Gamma}$, normalized to have
measure $1$ at each connected component of $\hat{\Gamma}$. Let $\nu$
be the law of the $2n$-tuple $X_{1},\ldots,X_{2n}$, $X_{2k}=\gamma_{k}+
\gamma_{k}^{-1}$,
$X_{2k-1}=-i(\gamma_{k}-\gamma_{k}^{-1})$. Then \[
\delta_{c}(\nu)=\dim_{\mathbb{C}}H^{1}(\Gamma;\mathbb{C}.)\]
\end{thm}
\begin{proof}
Let $\Gamma=\Gamma_{1}\oplus\Gamma_{2}$, where $\Gamma_{1}$ is a
finite group of order $l$ and $\Gamma_{2}$ is a free abelian group
on $p$ generators. Then $\hat{\Gamma}=\Gamma_{1}\times\mathbb{T}^{p}$,
where $\mathbb{T}$ denotes the unit circle in the complex plane.
Since $\mu$ is the Haar measure on $\hat{\Gamma}$, it is dimension
regular of dimension $p$. Hence $\delta_{c}(\nu)=lp$.

On the other hand, $H^{1}(\Gamma;\mathbb{C})=H^{1}(\hat{\Gamma};\mathbb{C}^{p})=\mathbb{C}^{p}$
and thus also has dimension $lp$. 
\end{proof}

\section{$\delta_{c}$ via Fisher information and a notion of Ricci curvature.}

\label{ricciclas} In this section, we relate $\delta_{c}$ with quantities
related with differential calculus. Let us remark, in the spirit of
Voiculescu \cite{Vo5}, that we can express $\delta_{c}$ via the
asymptotics of the associated Fisher information. To that end, recall
that for a probability measure $\mu(dx)=p(x)dx$ absolutely continuous
with respect to Lebesgue measure, the Fisher information is given
by \[
F(\mu)=\int(\partial_{x}\log p(x))^{2}p(x)dx.\]
 Note that if $P_s\mu=\mu_{\sqrt{s}}=\mu*p_{s}$ with $p_s=\nu_{\sqrt{s}
}$ the centered
Gaussian law with covariance $s$, since $\partial_s\frac{dP_s\mu}{dx}=
\frac{1}{2} (\frac{dP_s\mu}{dx})''$, 
$\partial_{s}H(P_{s}\mu)=-\frac{1}{2}F(P_{s}\mu)$
from which one sees that the entropy $H$ and the Fisher information
$F$ are related by \[
H(\mu)-H(\mu_{1})=\frac{1}{2}\int_{0}^{1}F(P_{s}\mu)ds.\]
 Taking $\mu=P_{t}\mu$ gives, since $H(\mu_{1})$  is always
bounded, that \begin{equation}
\delta_{c}(\mu)=1-\liminf_{t\rightarrow0}\frac{\int_{t}^{1}F(P_{s}\mu)ds}{
2|\log t^{\frac{1}{2}}
|}=1-\liminf_{t\rightarrow0}\frac{\int_{t}^{1}F(P_{s}\mu)ds}{
|\log t
|}
.\label{deltacfisher}\end{equation}
Observe that if $p_{s}$ is the density of $P_{s}\mu$ \[
\partial_{x}\log p_{s}(x)=\frac{1}{\sqrt{s}}E[g|X+\sqrt{s}g]\]
 when $g$ is a standard Gaussian variable independent from $X$ with
law $\mu$. This shows  by Cauchy-Schwartz inequality that \begin{equation}
0\le F(P_s\mu)\le\frac{1}{s}\label{boundfisher}\end{equation}
 and so proves again that $0\le\delta_{c}(\mu)\le1$. Moreover, \eqref{deltacfisher}
already reveals that $\delta_{c}(\mu)$ is related with the behaviour
of the Fisher information of $P_t\mu$ for small $t$ and in fact,
with the way that $P_t\mu$ approaches $\mu$ as $t$ goes to zero.
Let us give some heuristics by assuming that we have the stronger
statement that \[
F(P_t\mu)\approx_{t\rightarrow0}\frac{1-\delta_{c}(\mu)}{t}(1+o(1))\]
 and show that this entails that the convergence of $P_t\mu$ towards
$\mu$ is at least of the order $\sqrt{(1-\delta_{c}(\mu))t}$. In fact,
Fisher's information can be equivalently defined by \[
F(P_t\mu):=2\sup_{f}\{P_t\mu(\Delta f)-\frac{1}{2}P_t\mu((f')^{2})\}=\sup_{f}\frac{(P_t\mu(\Delta f))^{2}}{P_t\mu((f')^{2})}\]
 where the supremum is taken over all twice differentiable functions
$f$ (and is achieved here at $\log p_{t}$). Consequently, we find
that for all twice differentiable function $f$, \[
(P_t\mu(\Delta f))^{2}\le F(P_t\mu)\| f'\|_{\infty}^{2}.\]
 As a consequence, \begin{eqnarray*}
|P_t\mu(f)-\mu(f)| & \le & \int_{0}^{t}|\partial_{s}P_s\mu(f)|ds\\
 & = &\frac{1}{2} \int_{0}^{t}|P_s\mu(\Delta f)|ds\\
 & \le & \frac{1}{2}\| f'\|_{\infty}\int_{0}^{t}\sqrt{\frac{1-\delta_{c}(\mu)}{s}(1+o(1))}ds\\
 & \le &\| f'\|_{\infty}\sqrt{(1-\delta_{c}(\mu))t}(1+o(1)).\end{eqnarray*}
 Extending this inequality to all Lispchitz functions gives a bound
on the Duddley distance between $P_t\mu$ and $\mu$; \[
d(P_t\mu,\mu):=\sup_{f\textrm{ Lipschitz with norm}\leq1}|P_t\mu(f)-\mu(f)|\le
\sqrt{(1-\delta_{c}(\mu))t}(1+o(1)).\]
 We believe that the relation between the short time asymptotics of
$P_t\mu$ and $\delta_{c}$ should be deeper that this result even
though we could not prove it here. However, we shall prove here another
definition for $\delta_{c}$ which is closely related with Bochner's
inequality, a classical tool to estimate the short time asymptotics
of the heat kernel in a compact Riemaniann manifold. We shall restrict
ourselves here to measures on $\R$ but could as well consider measures
on a compact Riemaniann manifold with Ricci curvature bounded below
(eventually by a negative real number). To make this generalization
more transparent, we denote $\Delta$ the Laplace Baltrami operator
on $\RR$ (i.e the second spatial derivative). We let $\Gamma$ be
the carr\'{e} du champ given by\[
\G(f,g)=\frac{1}{2}\left( \D(fg)-f\D g-g\D f\right),\]
 and $\G_{2}$ be the carr\'{e} du champ it\'{e}r\'{e} \[
\G_{2}(f,f)=\frac{1}{2}\left(\D\G(f,f)-2\G(f,\D f)\right).\]
 In the case where $M=\RR$, we simply have \[
\G(f,f)=(f')^{2},\quad\G_{2}(f,f)=(f'')^{2}.\]
 Note that in the case of a connected Riemanian manifold with metric
$g$, Laplace Baltrami operator $\Delta$ and gradient $\nabla$,
the same definitions hold and give \[
\G(f,f)=g(\nabla f,\nabla f),\Gamma_{2}(f,f)=(\mbox{Hess}f,\mbox{Hess}f)_{g}+\mbox{Ric}(\nabla f,\nabla f)\]
 with $Ric$ the Ricci tensor. Bochner's (or curvature-dimension)
inequality $CD(n,K)$ states that \[
\Gamma_{2}(f,f)(x)\ge\frac{1}{n}(\D f)^{2}(x)-K\Gamma(f,f)(x)\]
 for all smooth function $f$ and at all points $x$ of the manifold.
$n$ corresponds to the dimension of the manifold whereas the best
constant $-K$ corresponds to the smallest eigenvalue of the Ricci
tensor. It is well known (see Bakry and Ledoux \cite{BL}, Bakry and
Qian \cite{BQ} etc) that the coefficient $n$ governs the short time
scaling of the heat kernel (as $t^{-\frac{n}{2}}$). Here $n\ge0$
and $K$ is a real number which we will assume finite for a while.
In the real one dimensional case, we clearly have $K=0$ and $n=1$,
but the constant $n$ of course is universal and does not depend on
any measure. We next define the measure-dependent Bochner inequality
as follows.

%\begin{notation}
%We write $P_\e\mu=P_{\e}*\mu$, where $P_{\e}$ is the Gaussian measure
%of variance $\e$. 
%\end{notation}
\begin{defn}
We say that a probability measure $\mu$ on $\RR$ satisfies Bochner's
inequality with constants ${\rm CD_{m}}(K,n)$ if there exists $\d>0$
so that for all $0\le\e'\le\d$, all smooth functions $f$, \[
P_{\e'}\mu(\Gamma_{2}(f,f))\ge\frac{1}{n}[P_{\e'}\mu(\D f)]^{2}-K(\e',n)P_{\e'}\mu(\Gamma(f,f)).\]

\end{defn}
In the sequel, it will appear that interesting cases appear when the
constant $K(n,\e')$ may blow up with $\e'$, reason why $K$ will
be later some \textit{non negative} arbitrary function. $n$ is some
positive real number. 

\textbf{Remark.} Note here that assuming that Bochner's inequality
is true in expectation would lead to the stronger definition

\[
P_{\e'}*\mu(\Gamma_{2}(f,f))\ge\frac{1}{n}P_{\e'}*\mu[(\D f)^{2}]-K(\e',n)P_{\e'}*\mu(\Gamma(f,f)).\]
 However, the idea is that what we want is that the points belonging
to the microstates \[
\Gamma_{\d,\mu}:=\{ x_{1},\cdots,x_{N}:d(\frac{1}{N}\sum\d_{x_{i}},\mu)<\d\}\]
 approximately satisfy Bochner's inequality when $N$ goes to infinity
and $\e$ goes to zero. Applying the classical Bochner's inequality
to functions of the form $F(x_{1},\cdots,x_{N})=N^{-1}\sum f(x_{i}+\e g_{i})$
for independent standard Gaussian variables $(g_{1},\cdots,g_{N})$,
$\e>0$ and letting $N$ go to infinity gives our actual definition
of measure-dependent Bochner's inequality. Hence, roughly speaking,
$(n,-K(\e,n))$ represent the dimension and the smallest eigenvalue
of the Ricci tensor of a manifold where the entries $(x_{1}+\sqrt{\e}g_{1},\cdots,x_{N}+\sqrt{\e}g_{N})$
live when the $(x_{1},\cdots,x_{N})$ belong to $\Gamma_{\d,\mu}$,
for $\d$ arbitrarily small.

Based on measure-dependent Bochner's inequalities we shall now define
a new entropy dimension

\begin{defn}
Let $\mu$ be a probability measure on $\RR$. We define the $CD$-
dimension as \[
\d^{\square}(\mu):=1-\inf_{\mu\textrm{ satisfies }{\rm CD_{m}}(n,K)}(\liminf_{\e\ra0}
\frac{\int_{\e}^{1}K(y,n)dy}{\log\e^{-1}}+1)n.\]
Above, the infimum is taken over all couple $(n,K(.,n))$ 
such that $\mu$ satisfies ${\rm CD_{m}}(n,K)$.
\end{defn}
We now prove that $\d^{\square}$ equals $\d_{c}$. We first prove
that

\begin{lem}
\label{indelta1} For any probability measure $\mu$ on $\R^{d}$,
\[
\d^{\square}(\mu)\le\d_{c}(\mu).\]
 
\end{lem}
\begin{proof}
Note that for $d=1$, $(\Delta f)^{2}=\Gamma_{2}(f,f)$ but that the
following argument will generalize to dimension $d$ by Cauchy-Schwartz
inequality which gives $d\Gamma_{2}(f,f)\ge(\Delta f)^{2}$. Integrating
with respect to $\mu$ implies that for all $\e\ge0$ \[
[P_\e\mu(\Delta f)]^{2}\le
P_\e\mu[(\Delta f)^{2}]\le
P_\e\mu[\Gamma_{2}(f,f)].\]
 On the other hand, with $p_{\e}$ the density of $P_\e\mu$ with
respect to Lebesgue measure,

\begin{eqnarray*}
[P_\e\mu(\Delta f)]^{2} & = & \left(P_\e\mu[f'(\log p_{\e})']\right)^{2}\\
 & \le & P_\e\mu[(f')^{2}]P_\e\mu[((\log p_{\e})')^{2}]\\
 & = & P_\e\mu[\Gamma_{1}(f,f)]F(P_\e\mu)\end{eqnarray*}

Therefore, for all $\alpha\in[0,1]$, we have \[
[P_\e\mu(\Delta f)]^{2}\le\alpha
P_\e\mu[\Gamma_{2}(f,f)]+(1-\alpha)F(P_\e\mu)P_\e\mu[\Gamma_{1}(f,f)]\]
 and so $\mu$ satisfies $CD_{m}(n,K)$ with $n=\alpha$ and \[
K(\e,n)=n^{-1}(1-n)\ F(P_\e\mu )\]
 for all $\alpha\in[0,1]$. Then, \[
\liminf_{\e\ra0}(\log\e^{-1})^{-1}\int_{\e}^{1}K(y,n)dy\le(1-n)n^{-1}\liminf_{\e\ra0}
(\log\e^{-1})^{-1}\int_{\e}^{1}F(P_y\mu)dy,\]
 and so \[
\d^{\square}(\mu)\ge1-\inf_{n\le
d}[n+(1-n)\liminf_{\e\ra0}(\log\e^{-1})^{-1}\int_{\e}^{1}F(P_y\mu)dy]=\d_{c}(\mu)\]
 where we used $(\log\e^{-1})^{-1}\int_{\e}^{1}F(P_y\mu)dy\le d=1$
by \eqref{boundfisher} to say that the infimum is taken at $n=0$.
\end{proof}
\begin{prop}
\label{bochner1} If a probability measure $\mu$ on $\RR$ satisfies
${\rm CD_{m}}(K,n)$, then

\[
\liminf_{\e\ra0}(\log\e^{-1})^{-1}\int_{\e}^{1}F(P_y\mu)dy\le\liminf_{\e\ra0}[(
\log\e^{-1})^{-1}\int_{\e}^{1}K(y,n)dy+1]n.\]
 
\end{prop}
As an immediate corollary of Proposition \ref{bochner1} we have

\begin{thm}
For any probability measure $\mu$ on $\RR$, \[
\d^{\square}(\mu)=\d_{c}(\mu)\]

\end{thm}
Whereas it can be easily seen that the characteristic $(n,-K)$ of
a manifold are invariant by Lipschitz map (simply by taking local
quadratic functions), invariance is not so transparent for measure-dependent
Bochner's inequality and we could not prove interesting invariance
property of $\d^{\square}$. However, the above theorem and section
\ref{secinv}  show that $\d^{\square}$ is invariant under Lipschitz
maps.

\begin{proof}
Let us first put $P_\e\mu=P_{\e}*\mu$ with $\e>0$ and write \[
F(P_\e \mu)=2\sup_{f}\{ P_\e\mu(\D f)-\frac{1}{2}P_{\e}\mu(\Gamma(f,f))\}\]
 Now, let for $x\in[0,\d]$, $\phi(x)=P_{x}*P_\e\mu(\Gamma(P_{\d-x}f,P_{\d-x}f))$ with $P_\e f(x)=P_\e(f(x)dx)$ by definition.
Differentiating with respect to $x$, we find that \begin{eqnarray*}
\phi'(x) & = & P_{x}*P_\e\mu(\Gamma_{2}(P_{\d-x}f,P_{\d-x}f))\\
 & \ge & \frac{1}{n}[P_{x}*P_\e\mu(\D P_{\d-x}f)]^{2}-K(x+\e,n)P_{x}*P_\e
\mu(\Gamma(P_{\d-x}f,P_{\d-x}f))\\
 & = & \frac{1}{n}((P_\e\mu\D P_{\d}f)^{2})-K(x+\e,n)\phi(x)\end{eqnarray*}
 where we used the fact that $P_{x}$ is a semigroup which commutes
with the Laplacian. Also, we have used our measure-dependent Bochner's
inequality with $f\ra P_{\d-x}f$ and $\e'=x+\e$. We set $L(x)=e^{\int_{x}^{1}K(y,n)dy}$.
Integrating $x\in[0,\d]$, we deduce that \begin{eqnarray}
P_{\d}*P_\e\mu(\Gamma(f,f)) & \ge & P_\e\mu
(\Gamma(P_{\d}f,P_{\d}f))\frac{L(\e+\d)}{L(\e)}+
\frac{1}{n}P_\e\mu((\D P_{\d}f))^{2}\int_{0}^{\d}\frac{L(\e+\d)}{L(\e+x)}dx\label{ineq}\end{eqnarray}
 We thus obtain that for all $a\in[0,1]$, \begin{eqnarray}
F(P_{\e+\d}) & \le & 2\sup_{f}\{ aP_\e \mu(\Delta P_{\d}f)-\frac{1}{2}P_\e
\mu(\Gamma(P_{\d}f,P_{\d}f))\frac{L(\e+\d)}{L(\e)}\nonumber \\
 &  & +(1-a)P_\e \mu(\Delta P_{\d}*f)-\frac{1}{2n}\int_{0}^{\d}\frac{L(\e+\d)}{L(\e+x)}dx(P_\e\mu(\D P_{\d}f))^{2}\}\nonumber \\
 & \le & a^{2}\frac{L(\e)}{L(\e+\d)}F(P_\e \mu)+(1-a)^{2}\frac{n}{\int_{0}^{\d}\frac{L(\e+\d)}{L(\e+x)}dx.}\label{bo1}\end{eqnarray}
 The optimum with respect to $a$ is taken at \[
a=\frac{n}{\frac{L(\e)}{L(\e+\d)}\int_{0}^{\d}\frac{L(\e+\d)}{L(\e+x)}dx
F(P_{\e}\mu)+n}.\]
 We conclude \begin{eqnarray}
F(P_{\e+\d}\mu) & \le & \frac{n\frac{L(\e)}{L(\e+\d)}F(P_\e\mu)}{\int_{0}^{\d}\frac{L(\e)}{L(\e+x)}dxF(P_\e\mu)+n}\label{bo2}\\
 & = & n\partial_{\d}[\log(\int_{0}^{\d}\frac{L(\e)}{L(\e+x)}dxF(P_\e\mu)+n)].\nonumber \end{eqnarray}
 Integrating with respect to $\d\in[0,1-\e]$ thus gives \begin{eqnarray*}
n^{-1}\int_{\e}^{1}F(P_x\mu)dx & \le & \log(n^{-1}\int_{\e}^{1}\frac{L(\e)}{L(x)}dxF(P_\e\mu)+1)\\
 & \le & \log(n^{-1}\e^{-1}\int_{\e}^{1}\frac{L(\e)}{L(x)}dx+1)\end{eqnarray*}
 where we used again $\e F(P_\e\mu)\le 1$ by \eqref{boundfisher}.
Consequently \begin{eqnarray*}
n^{-1}\liminf_{\e\ra0}(\log\e^{-1})^{-1}\int_{\e}^{1}F(P_x\mu)dx & \le & 
\liminf_{\e\ra0}(\log\e^{-1})^{-1}\log(n^{-1}\e^{-1}d\int_{\e}^{1}\frac{L(\e)}{L(x)}dx+1)\\
 & = &
\liminf_{\e\ra0}(\log\e^{-1})^{-1}\log(\e^{-1}\int_{\e}^{1}\frac{L(\e)}{L(x)}dx)\end{eqnarray*}
 Now, \[
\int_{\e}^{1}\frac{L(\e)}{L(x)}dx\le e^{\int_{\e}^{1}K(y,n)dy}\]
 and so we arrive at \begin{eqnarray}
n^{-1}\liminf_{\e\ra0}(\log\e^{-1})^{-1}\int_{\e}^{1}F(P_x\mu)dx & \le & 1+
\liminf_{\e\ra0}(\log\e^{-1})^{-1}\int_{\e}^{1}K(y,n)dy\end{eqnarray}
 which is the desired inequality.
\end{proof}
We finally give a lower bound of $\d^{\square}$ in the spirit of
\cite{Sh}. To do this, let us defined, for a $\Ca_{1}^{b}(\R,\R)$
function $g$, \[
F_{g}(\mu)=2\sup_{f}\{\mu(g\D f)-\frac{1}{2}\mu(\Gamma_{1}(f,f))\}\]

\begin{prop}
\label{proj} For any probability measure $\mu$ on $\RR^{d}$, \[
\d_{c}(\mu)=\d^{\square}(\mu)\ge1-\inf_{h\in\bar{\Fa}_{\mu}}\mu[(1-h)^{2}]\]
 with $\bar{\Fa}_{\mu}$ the set of continuous functions so that \[
\liminf_{\d\ra\infty}(\log\d^{-1})^{-1}\int_{\d}^{1}F_{h}(P_x\mu)dx=0.\]

\end{prop}
This lower bound has the advantage to give a more intuitive picture
of the dimension; for instance, if $\mu$ has a smooth density such
that the gradient of its logarithm is uniformly bounded, on a subset
$A$ of $M$, we take $h=1$ in some interior set $A^{s}$ of $A$,
$|h|\le1$ and $h=0$ outside $A$. It is easy to see that $F_{h}(\mu)<\infty$
and so $h\in\Fa_{\mu}$. Thus, we get \[
\d_c(\mu)=\d^{\square}(\mu)\ge\mu(A).\]
 Note however that such a lower bound is already contained in Theorem
\ref{thm:fractal}.

\begin{proof}
\noindent (of Proposition \ref{proj}). \textbf{}We take $h\in\Fa_{\mu}$.
We can assume without loss of generality that $\mu[(1-h)^{2}]\neq0$
since otherwise the bound is trivial ($h$ being equal to one almost
surely, and hence $F_{h}=F$ implying that $\d_{c}=\d^{\square}=d$).
We now write \begin{eqnarray*}
P_\e\mu(\D f) & = & P_\e\mu(h\D f)+P_\e\mu((1-h)\D f)\\
 & = & P_\e \mu(J_{h}f')+P_\e\mu((1-h)\D f)\end{eqnarray*}
 Now, \[
[P_\e\mu(J_{h}f')]^{2}\le F_{h}(P_\e\mu)P_\e\mu(\Gamma(f,f))\]
 whereas \begin{eqnarray*}
[P_\e\mu((1-h)\D f)]^{2} & \le & P_\e\mu
((1-h)^{2})P_\e\mu((\D f)^{2})\\
 & \le & P_\e\mu((1-h)^{2})P_\e\mu(\Gamma_{2}(f,f))\end{eqnarray*}
 Using that for all $\a>0$, for all $x,y\in\R$, $(x+y)^{2}\le(1+\a)x^{2}+(1+\a^{-1})y^{2}$
we thus derive the inequality\[
[P_\e\mu(\D f)]^{2}\le(1+\a)F_{h}(P_\e\mu)P_\e\mu(\Gamma(f,f))
+(1+\a^{-1})P_\e\mu((1-h)^{2})P_\e\mu(\Gamma_{2}(f,f))\]
that is the $CD_{m}(n,K)$ inequality with \[
n=n(\e)=(1+\a^{-1})P_\e\mu((1-h)^{2}),K(\e,n)=n^{-1}(1+\a)F_{h}(P_\e\mu)\]
 Since $h$ is continuous, $P_\e\mu((1-h)^{2})$ converges towards
$\mu((1-h)^{2})\neq0$ and since

\noindent
$\liminf(\log\e^{-1})^{-1}\int_{\e}^{1}F_{h}(P_x\mu)dx$
goes to zero , \[
\liminf_{\e\ra0}(\log\e^{-1})^{-1}\int_{\e}^{1}K(x,n)dx=0.\]
 Thus, $\d^{\square}(\mu)\ge1-\inf_{\alpha}(1+\a^{-1})\mu((1-h)^{2})=1-\mu((1-h)^{2})$
and optimizing over $h\in\bar{\Fa}_{\mu}$ yields the desired estimate.
\end{proof}

\section{Lipschitz invariance.}

\label{secinv}

Our main result is that $\delta_{c}$ is invariant under push-forwards
by bi-Lipschitz maps:

\begin{thm}
\label{thm:lip}Let $f:\mathbb{R}\to\mathbb{R}$ be bi-Lipschitz,
i.e., we assume that for some $m,M>0$ and all $x,y\in\mathbb{R}$,\[
m|x-y|\leq|f(x)-f(y)|\leq M|x-y|.\]
 Let $\eta=f^{*}\mu$ be the push-forward of $\mu$. Then $\delta_{c}(\mu)=\delta_{c}(\eta)$. 
\end{thm}
\begin{proof}
For any $y=f(x)$,\[
\eta[y-t/2,y+t/2]=\mu(f^{-1}[y-t/2,y+t/2])\geq\mu[x-t/(2M),x+t/(2M)].\]
 It follows that\begin{eqnarray*}
\int\log\frac{1}{t}\eta[y-t/2,y+t/2]d\eta(y) & \geq & \int\log\frac{1}{t}\mu[f^{-1}(y)-t/(2M),f^{-1}(y)+t/(2M)]d\eta(y)\\
 & = & \int\log\frac{1}{t}\mu[x-t/(2M),x+t/(2M)]d\mu(x)\\
 & = & \int\log\frac{1}{t/M}\mu[x-t/(2M),x+t/(2M)]d\mu(x)-\log M.\end{eqnarray*}
 Using Theorem \ref{thm:anynu} we conclude that\[
\delta_{c}(\eta)\leq\delta_{c}(\mu).\]
 Replacing $f$ by its inverse yields the reverse inequality. 
\end{proof}
It should be noted that one cannot expect much more invariance for
$\delta_{c}$ than is given by Theorem \ref{thm:lip}. Indeed, Cantor
sets in $\mathbb{R}$ can be made homeomorphic in a way that distorts
their fractal dimensions. \def\C{\mathbb{C}}

\section{Non-commutative Bochner's inequality}

\label{riccinc} In this last section, we generalize the notion of
measure-dependent Bochner's inequality of section \ref{ricciclas}.
To this end, we first define the appropriate notions of carr\'{e}
du champ and carr\'{e} du champ it\'{e}r\'{e}.

\subsection{Carr\'{e} du champ}

We recall first that the carr\'{e} du champ and the carr\'{e} du
champ it\'{e}r\'{e} in $\RR^{n}$ are given, for $f:\R^{n}\rightarrow\C$
by \[
\G(f,f)=\sum_{i=1}^{n}|\partial_{i}f|^{2},\quad\G_{2}(f,f)=\sum_{i,j=1}^{n}|\partial_{x_{i}}\partial_{x_{j}}f|^{2}\]
 In the case of $m$ Hermitian matrices $X_{N}$ with complex entries
$x_{ij}^{k}$, $1\le i\le j\le N$, $1\le k\le m$, \[
\D=2\sum_{k=1}^{m}\sum_{1\le i<j\le N}\partial_{x_{ij}^{k}}\partial_{\bar{x}_{ij}^{k}}+\sum_{k=1}^{m}\sum_{1\le i\le N}\partial_{x_{ii}^{k}}\partial_{x_{ii}^{k}}\]
 and so, if $f,g:\R^{2mN^{2}}\ra\mathbb{\C}$, we set \[
\G_{1}(f,g)=2\sum_{k=1}^{m}\sum_{i<j}\partial_{x_{ij}^{k}}f\partial_{\bar{x}_{ij}^{k}}\bar{g}+\sum_{k=1}^{m}\sum_{i<j}\partial_{x_{ii}^{k}}f\partial_{x_{ii}^{k}}\bar{g}\]
 and \[
\G_{2}(f,g)=\sum_{k,l=1}^{m}\sum_{ij}\sum_{ml}(\partial_{x_{ij}^{l}}\partial_{x_{ml}^{k}}f\partial_{\bar{x}_{ij}^{l}}\partial_{\bar{x}_{ml}^{k}}\bar{g}).\]
 Again, to define the notion of carr\'{e} du champ and carr\'{e}
du champ it\'{e}r\'{e} for tracial states, the idea is that if we
consider $f((x_{ml}^{k})_{1\le m\le l\le N}^{1\le k\le m}):=F(X)=\tr(P(X_{1},\cdots,X_{m}))$
when $\mun(Q):=N^{-1}\tr(Q(X_{1},\cdots,X_{m}))$ goes to $\tau(Q)$
for all polynomial $Q$ and some non-commutative law $\tau$. We denote
$*$ the involution \[
(zX_{i_{1}}\cdots X_{i_{k}})^{*}=\bar{z}X_{i_{k}}\cdots X_{i_{1}}\]
 for any $i_{l}\in\{1,\cdots,m\}$. Since $\overline{\tr(P)}=\tr(P^{*})$,
applying the above recipe we find,

\begin{eqnarray*}
\Gamma_{1}^{\mun}(P,Q) & =: & \sum_{k}\sum_{i,j}\partial_{x_{ij}^{k}}(\tr(P(X_{1},\cdots,X_{m})))\partial_{\bar{x}_{ij}^{k}}(\tr(Q^{*}(X_{1},\cdots,X_{m})))\\
 & = & N^{-1}\sum_{k}\sum_{i,j}[D_{k}P(X)]_{ij}[D_{k}Q(X)^{*}]_{ji}\\
 & = & \sum_{k}N^{-1}\tr(D_{k}P(X)(D_{k}Q(X))^{*})\\
 & \approx & \sum_{k}\tau(D_{k}P(X)(D_{k}Q(X))^{*}):=\Gamma_{1}^{\tau}(P,Q)\end{eqnarray*}
 where we have denoted by $D_{k}$ the cyclic derivative on polynomial,
given by$$D_{k}P=\sum_{P=P_{1}X_{k}P_{2}}P_{2}P_{1}$$ if $P$ is a
monomial (and extending by linearity to all polynomial then), and noticed,
as can be readily checked on monomials, that $(D_{k}P)^{*}=D_{k}P^{*}$.
Similarly,

\begin{eqnarray*}
\Gamma_{2}^{\mun}(P,Q) & =: & \sum_{k,l=1}^{m}\sum_{i,j}\sum_{pq}\partial_{X_{ij}^{l}}\partial_{X_{pq}^{k}}(\tr(P(X_{1},\cdots,X_{m})))\partial_{\bar{X}_{ij}^{l}}\partial_{\bar{X}_{pq}^{k}}(\tr(Q^{*}(X_{1},\cdots,X_{m})))\\
 & = & N^{-2}\sum_{k,l=1}^{m}[\partial_{l}\circ D_{k}P\sharp1_{pq}]_{ij}[\partial_{l}\circ D_{k}Q^{*}\sharp1_{qp}]_{ji}\\
 & \approx & \sum_{k,l=1}^{m}\tau\otimes\tau((\partial_{l}\circ D_{k}Q)^{*}\star\partial_{l}\circ D_{k}P):=\Gamma_{2}^{\tau}(P,Q)\end{eqnarray*}
 where $\partial_{k}$ denotes the non-commutative derivative with
respect to the variable $X_{k}$ ($\partial_{k}P=\sum_{P=P_{1}X_{k}P_{2}}P_{1}\otimes P_{2}$
for a monomial $P$), $A\otimes B\sharp C=ACB$, $(A\otimes B)^{*}=B^{*}\otimes A^{*}$
and $A\otimes B\star A'\otimes B'=BA'\otimes AB'$. $1_{kl}$ is the
matrix with zeroes except in $kl$. Hence, we define 

\begin{defn}
For any non-commutative law $\tau$ of $m$ self-adjoint variables,
we define its \textit{non-commutative carr\'{e} du champ}  to be
the bilinear function on $\cxm$ so that for any $P,Q\in\cxm$, \[
\Gamma_{1}^{\tau}(P,Q)=\sum_{i=1}^{m}\tau(D_{i}P(D_{i}Q)^{*})\]
 and its \textit{non-commutative carr\'{e} du champ it\'{e}r\'{e}}
to be the bilinear function on $\cxm$ so that for any $P,Q\in\cxm$,

\[
\Gamma_{2}^{\tau}(P,Q)=\sum_{k,l=1}^{m}\tau\otimes\tau((\partial_{l}\circ D_{k}Q)^{*}\star\partial_{l}\circ D_{k}P).\]
 We also denote in short \[
\Gamma_{i}^{\tau}(P,Q)=<P,Q>_{\tau,i}.\]
 Observe that the above notation makes sense since $\Gamma_{i}^{\tau}$
are positive bilinear forms. This is obvious for $\Gamma_{1}^{\tau}$.
For $\Gamma_{2}^{\tau}$, one needs to observe that if $\tau$ is
a tracial state, $P,Q\ra\tau\otimes\tau(P\star Q^{*})$ is non negative.
But if $P=\sum\a_{i}A_{i}\otimes B_{i}$, \[
\tau\otimes\tau(P\star P^{*})=\sum\a_{i}\bar{\a}_{j}\tau(A_{i}A_{j}^{*})\tau(B_{i}B_{j}*)\ge0\]
 since the matrices $(\tau(A_{i}A_{j}^{*}))_{i,j},(\tau(B_{i}B_{j}^{*}))_{i,j}$
are non-negative.
\end{defn}
Let us introduce the notation: 
\[\partial_{k}^{2}\equiv\frac{1}{2}\left(\partial_{k}\otimes1+1\otimes\partial_{k}\right)\circ\partial_{k},\]
\[
M(A\otimes B\otimes C)\equiv B\otimes AC\]
 and \[
{\mathbb{L}}_{\tau}:=\sum_{k}(\tau\otimes I)(M\circ\partial_{k}^{2}).\]
Then when the entry-wise Laplacian $\Delta=\sum\partial_{x_{ij}^{k}}\partial_{\bar{x}_{ij}^{k}}$
acts on $F(X_{ij}^{l})=f(X_{1},\cdots,X_{m})$, we
get 
 that \[
\Delta F={\mathbb{L}}_{\tau}f\]
 when the law of $X$ approximates $\tau$. If $F=N^{-1}\tr(P)$, we
get \[
\Delta F\approx\tau({\mathbb{L}}_{\tau}P).\]
 Note here that \[
\tau({\mathbb{L}}_{\tau}P)=\sum_{i=1}^{m}\tau\otimes\tau(\partial_{i}\circ D_{i}P)\]
 as can be readily checked by taking $P$ to be a monomial. Let $S=(S^{1},\cdots,S^{m})$
be a free Brownian motion, free with $X=(X^{1},\cdots,X^{m})$ with
law $\tau$, and $\phi$ a tracial state on a von Neumann algebra
containing $S$ and $X$. We then have \[
P(X+S_{t})=P(X)+\int_{0}^{t}{\mathbb{L}}_{\phi_{X+S_{s}}}(P)(X+S_{s})ds+\int_{0}^{t}\sum_{i=1}^{m}\partial_{i}P(X+S_{s})\sharp dS_{s}^{i}\]
 where the last term is a martingale. We denote $\tau_{t}$ the distribution
of $(X^{1}+S_{t}^{1},\cdots,X^{m}+S_{t}^{m})$.

\subsection{Non-commutative Bochner's inequality}

We recall that Bochner's inequality reads in the classical context
as \[
\Gamma_{2}(f,f)\ge\frac{1}{n}(\D f)^{2}-K\Gamma_{1}(f,f)\]
 for some fixed constants $n\ge0$, $K\in{\mathbb{R}}$. Remark that
$n$ is of the order of the dimension, so of order $N^{2}$ in the
context of matrices, so we let $\Na=n/N^{2}$ and apply this inequality
to $F=\tr(P)$ we get if $\mun_{X}\approx\tau$, as $N$ goes to infinity,
\[
<P,P>_{\tau,2}\ge\frac{1}{\Na}[\tau({\mathbb{L}}_{\tau}P)]^{2}-K<P,P>_{\tau,1}.\]
 Therefore,

\begin{defn}
We shall say that a non-commutative law $\tau$ satisfies a ${\rm CD_{m}}(\Ka,\Na)$
inequality iff for all $\e$ small enough,

\[
<P,P>_{\tau_{\e},2}\ge\frac{1}{\Na}[\tau_{\e}({\mathbb{L}}_{\tau_{\e}}P)]^{2}-\Ka(\Na,\e)<P,P>_{\tau_{\e},1}\]
for any polynomial function $P$.
 
\end{defn}
We can therefore define

\begin{defn}
\[
\d^{\square}(\tau)=m-\inf_{\tau\textrm{ satisfies }{\rm CD_{m}}(\Ka,\Na)}(\bar{\Ka}(\Na)+1)\Na\]
 where \[
\bar{\Ka}(\Na)=\liminf_{\e\ra0}(\log\e^{-1})^{-1}\int_{\e}^{1}\Ka(\Na,y)dy.\]

\end{defn}
We next want to compare this definition of a non-commutative dimension
with already existing entropy dimension. We recall that in the non-commutative
setting, Voiculescu \cite{Vo6} defined the following notion of Fisher
entropy and related entropy dimension. For a tracial state $\tau$,
we define its Fisher information by \begin{eqnarray*}
\Phi^{*}(\tau) & = & \sum_{i=1}^{m}\sup_{P\in\cxm}\{\tau\otimes\tau(\partial_{i}(P+P^{*}))-\tau(PP^{*})\}\\
 & = & \sup_{P\in\cxm^{m}}\{\sum_{i=1}^{m}\tau\otimes\tau(\partial_{i}(P_{i}+P_{i}^{*}))-\sum_{i=1}^{m}\tau(P_{i}P_{i}^{*})\}\end{eqnarray*}
 Then, as in \eqref{deltacfisher}, the microstates-free free entropy
dimension is given by

\begin{equation}
\delta^{*}(\mu)=m-\liminf_{t\rightarrow0}\frac{\int_{t}^{1}\Phi^{*}
(\tau_{s})ds}{|\log t|}.\label{deltancfisher}\end{equation}
 Here, we shall consider a variant of $\d^{*}$ based on the following
definition of Fisher information as found in \cite{CDG1}: \[
\bar{\Phi}^{*}(\tau)=\sup_{P\in\cxm}\{\sum_{i=1}^{m}
\tau\otimes\tau(\partial_{i}(D_{i}P+D_{i}P^{*}))-\sum_{i=1}^{m}
\tau(D_{i}PD_{i}P^{*})\}\]

and \[
\bar{\d}^{*}(\tau)=m-\liminf_{t\rightarrow0}
\frac{\int_{t}^{1}\bar{\Phi}^{*}(\tau_{s})ds}{|\log t|}.\]
 Observe that $\bar{\Phi}^{*}\le\Phi^{*}$ and so $\bar{\d}^{*}(\tau)\ge\d^{*}(\tau)$.
Equality is achieved if the conjugate variables belong to the cyclic
gradient space, which appears to be often (if not always) the case
(see Voiculescu \cite{Vo6} and Cabanal Duvillard-Guionnet \cite{CDG2}).
This is the case, in particular, if we are dealing with the law $\tau$
of a single variable (i.e., $m=1$).

In the sequel, we shall as well denote $(\Ja_{\tau}^{i})_{1\le i\le m}$
for the projection of the conjugate variable on the cyclic gradient
space, i.e \[
\tau\otimes\tau(\partial_{i}\circ D_{i}P)=\tau(\Ja_{\tau}^{i}D_{i}P)\]
 for all polynomials $P$.
We next prove 

\begin{prop}
\label{Bochnet} \[
\d^{\square}(\tau)=\bar{\d}^{*}(\tau).\]
 In particular, \[
\d^{\square}(\tau)\ge\bar{\d}^{*}(\tau)\ge\d(\tau)\]
 where $\delta(\tau)$ denotes the microstates entropy dimension.
\end{prop}
\begin{proof}
Let us first remark that by definition \[
\tau(\LL_{\tau}P)=\sum_{i=1}^{m}\tau\otimes\tau(\partial_{i}\circ D_{i}P)=
\sum_{i=1}^{m}\tau(\Ja_{\tau}^{i}D_{i}P)\]
 and therefore \[
|\tau(\LL_{\tau}P)|^{2}\le\bar{\Phi}(\mu)\Gamma_{1}^{\tau}(P,P).\]
 On the other hand \[
|\tau(\LL_{\tau}P)|^{2}\le m\sum_{i=1}^{m}|\tau\otimes\tau(\partial_{i}
\circ D_{i}P)|^{2}\]
 with \[
|\tau\otimes\tau(\partial_{i}\circ D_{i}P)|^{2}\le\tau\otimes\tau(\partial_{i}
\circ D_{i}P\star(\partial_{i}\circ D_{i}P)^{*})\]
 by Cauchy-Schwartz inequality, which holds because of the positivity
of the positive bilinear form $P,Q\ra\tau\otimes\tau(\partial_{i}\circ 
D_{i}P\star(\partial_{i}\circ D_{i}P)^{*})$.
Hence, for any $\a\in[0,1]$ \begin{eqnarray*}
|\tau(\LL_{\tau}P)|^{2} & \le & m\a\Gamma_{2}^{\tau}(P,P)+
(1-\alpha)\bar{\Phi}(\tau)\Gamma_{1}^{\tau}(P,P).\end{eqnarray*}
 This proves that Bochner's inequality is satisfied with $\Na=m\a$
and $\Ka(\Na,\e)=(1-\Na/m)\bar{\Phi}(\tau_{\e})\Na^{-1}$ from which
we get \[
m-\delta^{\square}(\tau)=\inf\{\Na(1+\bar{\Ka}(\Na))\}
\le\inf_{\Na\in[0,m]}\{\Na+(1-\Na/m)\liminf
\frac{\int_{\e}^{1}\bar{\Phi}^{*}(\tau_{s})ds}{|\log\e|}\}=m-\bar{\delta}^{*}(\tau)\]
 where we used that $\frac{\int_{\e}^{1}\bar{\Phi}^{*}(\tau_{s})ds}{|\log\e|}
\in[0,m]$
which holds since $\bar{\Phi}^{*}(\tau_{s})\le s^{-1}$.

For the other inequality, let $X$ be an $m$-tuple of random variables
having the law $\tau_{x+\e}$ obtained as free convolution of the
law $\tau$ with the semicircular law of variance $\e$. Let $0<x<\delta$
and let $S_{\delta-x}$ be an $m$-tuple of semicircular variables
of variance $\delta-x$, free from $X$. Denote by $\tau(\cdot|X)$
the conditional expectation onto the algebra generated by $X$. We
then introduce, in the spirit of the proof in the classical case,
the function \[
\phi(x)=\sum_{i=1}^{m}\tau_{x+\e}\left(|D_{i}\tau(P(X+S_{\delta-x})|X)|^{2}\right)\]
(note that $\tau(P(X+S_{\delta-x})|X)$ is a polynomial in $X$ and
hence is in the domain of $D_{i}$).

We have \begin{eqnarray}
\phi'(x) & = & \sum_{i=1}^{m}\tau_{x+\e}\left(\LL_{\tau_{x+\e}}|D_{i}\tau(P(X+S_{\delta-x})|X)|^{2}\right)\nonumber \\
 &  & -2\Re\tau_{x+\e}\left(D_{i}\tau(\LL_{\tau_{\d+\e}}P(X+S_{\delta-x})|X)(D_{i}\tau(P(X+S_{\delta-x})|X)^{*}\right)\label{eql}\end{eqnarray}
 where we used the fact that the law of $X+S_{\delta-x}$ under $\tau_{x+\e}$
is the law of $X+S_{\delta-x}+\bar{S}_{x+\e}$, with $\bar{S}$ a
free Brownian motion independent from $S,X$, which has the same law
$\tau_{\d+\e}$ of $X+S_{\delta+\e}$. Now, let us compute $\LL_{\tau_{x+\e}}(PQ)$
for polynomials $P,Q$. $\LL_{\tau_{x+\e}}$ is a second order differential
operator; it will either act on $P$, or $Q$, or both; \[
\LL_{\tau_{x+\e}}(PQ)=\LL_{\tau_{x+\e}}(P)Q+P\LL_{\tau_{x+\e}}(Q)+R(P,Q).\]
 To compute $R(P,Q)$ note that this contribution comes from \[
\D_{k}^{2}(PQ)-\D_{k}^{2}(P)\times1\otimes1\otimes Q-P\otimes1\otimes1\times\D_{k}^{2}(Q)=\partial_{k}P\bar{\star}\partial_{k}Q\]
 with $A\otimes B\bar{\star}A'\otimes B'=A\otimes BA'\otimes B'$.
Note that \[
M(A\otimes B\bar{\star}A'\otimes B')=BA'\otimes AB'=A\otimes B\star A'\otimes B'.\]
 Therefore \[
\sum_{i=1}^{m}\tau_{x+\e}\left(R(D_{i}\tau(P(X+S_{\delta-x})|X),D_{i}\tau(P(X+S_{\delta-x})|X))\right)\]
 \[
=\Gamma_{2}^{\tau_{x+\e}}(\tau(P(X+S_{\delta-x})|X),\tau(P(X+S_{\delta-x})|X))\]
 Finally, it is easy to see that \[
\LL_{\tau_{x+\e}}(D_{i}\tau(P(X+S_{\delta-x})|X))=D_{i}\tau(\LL_{\tau_{\d+\e}}P(X+S_{\delta-x})|X)\]
 so that we have proved according to \eqref{eql} that \begin{eqnarray}
\phi'(x) & = & \Gamma_{2}^{\tau_{x+\e}}(\tau(P(X+S_{\delta-x})|X))\nonumber \\
 & \ge & \frac{1}{\Na}[\tau_{x+\e}[\LL_{\tau_{x+\e}}(\tau(P(X+S_{\delta-x})|X)]^{2}-K\Gamma_{1}^{\tau_{x+\e}}(\tau(P(X+S_{\delta-x})|X))\end{eqnarray}
 We can now proceed exactly in the lines of the proof of Proposition
\ref{Bochnet} to conclude that $\bar{\Phi}^{*}(\tau_{\e})$ satisfies
the bound \begin{eqnarray}
\bar{\Phi}^{*}(\tau_{\e}) & \le & \frac{\Na\frac{L(\e)}{L(\e+\d)}\bar{\Phi}^{*}(\tau_{\e})}{\int_{0}^{\d}\frac{L(\e)}{L(\e+x)}dx\bar{\Phi}^{*}(\tau_{\e})+\Na}\label{bo22}\end{eqnarray}
 with $L(y)=e^{\int_{y}^{1}\Ka(x,\Na)dx}$ as before. The rest of
the proof is exactly as in the classical case.
\end{proof}
\begin{cor}
If $\tau$ is the law of a single variable (i.e., $m=1$) then\[
\d^{\square}(\tau)=\bar{\d}^{*}(\tau)=\d(\tau)=1-\tau\otimes\tau(\chi_{\Delta})\]
where $\chi_{\Delta}$ is the characteristic function of the diagonal
$\Delta\subset\mathbb{R}^{2}$ and we identify $\tau$ with a measure
on $\mathbb{R}$.
\end{cor}
\begin{prop}
Let $X=(X_{1},\ldots,X_{m})$ have the given law $\tau$ , $M=W^{*}(X_{1},\ldots,X_{m})$
and let $G=(G_{ij})\in M_{m\times m}(L^{2}(M\bar{\otimes}M^{o}))$
be a fixed matrix. Let $\bar{\Phi}_{G}$ be the Fisher information
defined by \[
\bar{\Phi}_{G}=\sup_{P\in\cxm}\{\sum_{i=1}^{m}\tau\otimes\tau(\partial_{i}^{G}(D_{i}P+D_{i}P^{*}))-\sum_{i=1}^{m}\tau(D_{i}PD_{i}P^{*})\}\]
where $\partial_{i}^{G}(X_{j})=G_{ij}$. Then\[
\bar{\delta}^{*}(\tau)=\d^{\square}(\tau)\ge m(1-\inf_{G\in\Fa_{\tau}}\tau(1-G)^{2})\]
 with $\Fa_{\tau}$ the set of $G\in M_{m\times m}(L^{2}(M\bar{\otimes}M^{o}))$
so that $(\log\e^{-1})^{-1}\int_{\e}^{1}dt\bar{\Phi}_{G}^{*}(\tau_{t})$
goes to zero. 
\end{prop}
The proof is exactly the same as the previous one except that the
use of Bochner inequality is simply replaced by the fact that any
measure satisfies ${\rm CD_{m}}(m,0)$ as we have seen in the proof
of the previous theorem.

{\bf Acknowledgment} Alice Guionnet wishes to
thank M. Ledoux  for motivating discussions  on part of this article.

\end{document}